\documentclass{svproc}
\usepackage{url}
\usepackage{amssymb}
\usepackage{enumerate}
\usepackage{amsmath}
\usepackage{mathrsfs}
\usepackage{xcolor}

\newcommand{\E}{\tilde{\mathbb{E}}}

\newcommand{\Dm}{\mathcal{U}(\mathfrak{D}_{n-1})}

\newtheorem{conjecture-non}{Conjecture}

\begin{document}
\mainmatter              
\title{The Racah algebra and $\mathfrak{sl}_n$}
\titlerunning{The Racah algebra and $\mathfrak{sl}_n$}  
%
\author{Hendrik De Bie\inst{1} \and Luc Vinet\inst{2}
 \and
Wouter van de Vijver\inst{1} }
\authorrunning{Hendrik De Bie et al.} 
%
%
\institute{Department of Electronics and Information Systems, Faculty of Engineering and Architecture, Ghent University, Krijgslaan 281, Building S8, 9000 Gent, Belgium\\
\email{hendrik.debie@ugent.be}, \email{wouter.vandevijver@ugent.be},\\
\and
Centre de Recherches Math\'ematiques, Universit\'e de Montr\'eal, P.O. Box 6128, Centre-ville Station, Montr\'eal, QC H3C 3J7, Canada\\ \email{vinet@crm.umontreal.ca}}

\maketitle              

\begin{abstract}

We conjecture the existence of an embedding of the Racah algebra into the universal enveloping algebra of $\mathfrak{sl}_n$. Evidence of this conjecture is offered by realizing both algebras using differential operators and giving an embedding in this realization.
\end{abstract}
\section{Introduction}

The Racah algebra synthecizes the properties of the Racah polynomials \cite{Genest&Vinet&Zhedanov-2014-3,Zhedanov-1991}, which are the most complicated univariate discrete orthogonal polynomials in the Askey scheme \cite{Koek}. 

Multivariate Racah polynomials were introduced by Tratnik in \cite{Trat}. These polynomials also have a solid algebraic underpinning, as was recently established in \cite{wouter} using the higher rank Racah algebra. This higher rank Racah algebra was initially introduced in \cite{I,I2} in the context of superintegrability and later in \cite{racah} as a subalgebra of intermediate Casimir elements in the $n$-fold tensor product of $\mathfrak{su}(1,1)$.

Although the initial motivation to introduce the (higher rank) Racah algebra was to establish a connection with the multivariate Racah polynomials, the algebra has now become an independent object of study. In particular its relation with other algebraic structures is part of ongoing investigations. We refer the reader to \cite{crampe2} for connections with Brauer and Temperley-Lieb algebras, and to \cite{gab} for connections with Howe duality.

The present paper aims to provide evidence for the following conjecture:
\begin{conjecture-non}\label{Conjecture}
There exists an embedding of the higher rank Racah algebra into the universal enveloping algebra of the Lie algebra $\mathfrak{sl}_n$.
\end{conjecture-non}

Indeed, we construct this embedding for a differential operator realization of $\mathcal{R}_n$ (recently introduced in \cite{barg}) in the enveloping algebra of a differential operator realization of $\mathfrak{sl}_n$. Note that the embedding in the rank one case was already constructed in \cite{equit}.

\section{Definition of the higher rank Racah algebra}

The algebra $\mathfrak{su}(1,1)$ is generated by three elements $A_\pm$ and $A_0$ with following relations:
\begin{equation*}\label{su11}
[A_-,A_+]=2A_0, \qquad [A_0, A_\pm]=\pm A_\pm.
\end{equation*}
Its universal enveloping algebra $\mathcal{U}(\mathfrak{su}(1,1))$ contains the Casimir element of $\mathfrak{su}(1,1)$:
\begin{equation*}\label{Casimir}
	C:=A_0^2-A_0-A_+A_-.
\end{equation*}
We define the following elements of $\mathcal{U}(\mathfrak{su}(1,1))^{\otimes n}$ for $1\leq k  \leq n$
\begin{align*}
	A_{0,k}&:= 1^{\otimes (k-1)} \otimes A_0 \otimes  1^{\otimes (n-k)},\quad \\
	A_{\pm,k}&:=1^{\otimes (k-1)} \otimes A_\pm \otimes  1^{\otimes (n-k)}.
\end{align*}
For any non-empty subset  $K \subset [n]:=\{1,\ldots, n \}$ we define similarly
\begin{equation*}
	A_{0,K}:=\sum_{k\in K} A_{0,k}, \quad A_{\pm,K}:=\sum_{k\in K} A_{\pm,k}.
\end{equation*}
The three operators $A_{0,K}$ and $A_{\pm,K}$ generate an algebra isomorphic to $\mathfrak{su}(1,1)$. Its Casimir is given by.
\begin{align*}
\label{Casimir-Upper}
	C_K:=A_{0,K}^2-A_{0,K}-A_{+,K}A_{-,K}.
\end{align*}
These operators generate the higher rank Racah algebra.
\begin{definition}
The higher rank Racah algebra $\mathcal{R}_n$ is the subalgebra of  $\mathcal{U}(\mathfrak{su}(1,1))^{\otimes n}$ generated by the set of operators
\[
	\{ C_{A}\, |\, A \subset \{1, \dots , n\} \text{ and } A\neq \emptyset \}
\]
\end{definition}
A full account of $\mathcal{R}_n$ is presented in \cite{racah}. We mention one fact here: the operators $C_A$ are not linearly independent. By formula $(17)$ in \cite{racah} we have:
\begin{equation*}
\label{Relation}
C_A=\sum_{\left\{i,j\right\}\subset A} C_{ij}-\left(|A|-2\right)\sum_{i \in A} C_i .
\end{equation*}
Hence, if one wants to present realizations of $\mathcal{R}_n$, it suffices to give expressions for the operators $C_{ij}$ and $C_i$. In section \ref{realization} we will present the higher rank Racah algebra as given in \cite{barg} this way.

\section{Realizing $\mathfrak{sl}_n$ in $n-1$ variables}
Let $\mathfrak{sl}_n(\mathbb{R})$  be the algebra of $n\times n$ matrices whose trace equals zero and with the commutator as Lie bracket. Let $E_{ij}$ be the matrix whose entries are equal to $0$ except for the entry on the $i$th row and $j$th column which equals $1$. Then the Lie algebra $\mathfrak{sl}_n(\mathbb{R})$ is generated by the set
\[ \{ E_{ij} | 1 \leq i, j \leq n \text{ and } i\neq j \} \cup \{ E_{ii}-E_{nn} | 1 \leq i \leq n-1\}. \]

Let $u_i$, $i \in \{ 1 \dots n-1\}$ be real variables. We introduce the differential operators:
\begin{align*}
  T_{ij}&:=-u_j\partial_i \qquad  i\neq j \text{ and } i,j< n \\
  T_{in}&:=-\partial_i \qquad \qquad i< n \\
  T_{nj}&:=u_j\E \qquad \qquad j< n \\
  \tilde{T}_d&:=-u_d\partial_d-\E  \qquad d <n
 \end{align*}
where the operator $\E$ is defined as
\[
\E:=\sum_{i=1}^{n-1} u_i\partial_i-k.
\]
Using again the commutator as Lie bracket we denote by $\mathfrak{D}_n$ the Lie algebra spanned by all the $\tilde{T}_d$ and $T_{ij}$.
The real number $k$ is a deformation parameter that leaves the algebra relations invariant. 
One observes that $\mathfrak{sl}_n(\mathbb{R})$ and $\mathfrak{D}_n$ are isomorphic. The isomorphism $\sigma$  is given by
\begin{align*}
\sigma(E_{ij})&=T_{ij} \\
\sigma(E_{dd}-E_{nn})&=\tilde{T}_d.
\end{align*}
Note that this isomorphism does not extend to their universal enveloping algebras.

\subsection{Some operators in $\mathcal{U}(\mathfrak{D}_n)$}
We introduce a number of operators we need later on. Observe that the operator $\E$ is in the universal enveloping algebra $\mathcal{U}(\mathfrak{D}_n)$ but not in $\mathfrak{D}_n$ because $\mathfrak{D}_n$ lacks the identity:
\[
	 \E=\frac{-1}{n}\left(k+\sum_{d=1}^{n-1} \tilde{T}_d\right).
 \]
Let $u_B:=\sum_{k \in B} u_k$. We have
\[
	u_B\E:=\sum_{j \in B} T_{nj}.
\]
We will also express $u_B\partial_\alpha$ in function of the generators:
\begin{equation}
\label{vgl}
	u_B\partial_\alpha:=-\delta_{\alpha B}(\tilde{T}_\alpha+\E)+\sum_{j \in B \backslash \alpha} -T_{\alpha j},
\end{equation}
where we introduced a new symbol standing for:
\[
	\delta_{\alpha B}:=	\begin{cases} 	0, \text{ if } \alpha \notin B \\
								1, \text{ if } \alpha \in B.
					\end{cases}	
\]
It is then easy to check the following Lemma.
\begin{lemma}\label{lemma1} The following holds
\begin{align*}
	[ u_B\E,\partial_{\alpha}]&=-u_B\partial_{\alpha}-\delta_{\alpha B}\E \\
	[ u_A, \partial_{\alpha}]&=-\delta_{\alpha A}.
\end{align*}
\end{lemma}

\section{Realization of $\mathcal{R}_{n}$ in $n-2$ variables}\label{realization}
In \cite{barg} an explicit differential operator realization of $\mathcal{R}_n$ was given in Theorem $5$. We repeat this theorem here.
\begin{theorem}
The space $\Pi_k^{n-2}$ of all polynomials of degree $k$ in $n-2$ variables carries a realization of the rank $n-2$ Racah algebra $\mathcal{R}_n$.
This realization is given explicitly by
\[
\widetilde{C_{i}}=\nu_i(\nu_i-1), \qquad i \in [n]
\]
and, for $i, j \in \{3, \ldots, n \}$,
\begin{align*}
\widetilde{C_{12}}&=- \left(k-1-\sum_{\ell=1}^{n-2} u_{\ell}\partial_{u_\ell} \right)  \left(-k-\partial_{u_1}+\sum_{\ell=1}^{n-2} u_{\ell}\partial_{u_\ell} \right) + 2 \nu_2 \left(k-\sum_{\ell=1}^{n-2} u_{\ell}\partial_{u_\ell} \right) \\
& \qquad - 2 \nu_1\left(-k-\partial_{u_1}+\sum_{\ell=1}^{n-2} u_{\ell}\partial_{u_\ell} \right) + (\nu_1+\nu_2)(\nu_1+\nu_2-1)\\
\widetilde{C_{1j}}&=- \left(1 - \sum_{\ell=1}^{j-2} u_{\ell} \right)^2 \left(k-1-\sum_{\ell=1}^{n-2} u_{\ell}\partial_{u_\ell} \right)   \left( \partial_{u_{j-2}}- \partial_{u_{j-1}} \right)\\
& \qquad + 2 \nu_j \left(1 - \sum_{\ell=1}^{j-2} u_{\ell} \right)\left(k-\sum_{\ell=1}^{n-2} u_{\ell}\partial_{u_\ell} \right) 
 - 2 \nu_1  \left(1 - \sum_{\ell=1}^{j-2} u_{\ell} \right)  \left( \partial_{u_{j-2}}- \partial_{u_{j-1}} \right)\\ & \qquad  + (\nu_1+\nu_j)(\nu_1+\nu_j-1)\\
\widetilde{C_{2j}}&= -\left( \sum_{\ell=1}^{j-2} u_{\ell} \right)^2 \left(1-k-\partial_{u_1}+\sum_{\ell=1}^{n-2} u_{\ell}\partial_{u_\ell} \right)   \left( \partial_{u_{j-2}}- \partial_{u_{j-1}} \right)\\
& \qquad + 2 \nu_j \left( \sum_{\ell=1}^{j-2} u_{\ell} \right) \left(k+\partial_{u_1}-\sum_{\ell=1}^{n-2} u_{\ell}\partial_{u_\ell} \right)  + 2 \nu_2 \left( \sum_{\ell=1}^{j-2} u_{\ell} \right)\left( \partial_{u_{j-2}}- \partial_{u_{j-1}} \right)   \\
& \qquad + (\nu_2+\nu_j)(\nu_2+\nu_j-1)\\
\widetilde{C_{ij}}&= - \left( \sum_{\ell=j-1}^{i-2} u_{\ell} \right)^2 \left( \partial_{u_{i-2}}- \partial_{u_{i-1}} \right) \left( \partial_{u_{j-2}}- \partial_{u_{j-1}} \right)\\
 & \qquad+ 2 \nu_j  \left( \sum_{\ell=j-1}^{i-2} u_{\ell} \right)  \left( \partial_{u_{i-2}}- \partial_{u_{i-1}} \right) - 2 \nu_i  \left( \sum_{\ell=j-1}^{i-2} u_{\ell} \right)  \left( \partial_{u_{j-2}}- \partial_{u_{j-1}} \right)\\
& \qquad + (\nu_i+\nu_j)(\nu_i+\nu_j-1)
\end{align*}
where we assume $i >j$ and with $u_{n-1}=0$ whenever it appears.
\end{theorem}

We want to express these operators as elements in $\mathcal{U}(\mathfrak{D}_{n-1})$. 

\subsection{The differential embedding}
To show that each generator of $\mathcal{R}_{n}$ is in $\mathcal{U}(\mathfrak{D}_{n-1})$, we wil express each generator in function of $\partial_\alpha$, $\E$, $u_B\partial_{\alpha}$ and $u_B\E$.
Let us start with the operator $\widetilde{C_{12}}$. We can express this operator as follows:
\begin{align*}
 \widetilde{C_{12}}&=-\left(-\E-1\right)\left(-\partial_1+\E\right)+2\nu_2\left(-\E\right)-2\nu_1\left(-\partial_1+\E\right)\\
& \quad+(\nu_1+\nu_2)(\nu_1+\nu_2-1)
\end{align*}
As $\partial_1=-T_{1n-1}$ and $\E$ are in $\mathcal{U}(\mathfrak{D}_{n-1})$ so is $\widetilde{C_{12}}$.

Consider the first term of the operator $\widetilde{C_{1j}}$:
\begin{align*}
	&-\left(1-u_{[j-2]}\right)^2\left(-1-\E\right)\left(\partial_{j-2}-\partial_{j-1}\right) \\
	&=\left(1-u_{[j-2]}\right)^2\left(\partial_{j-2}-\partial_{j-1}\right)\E \\
	&=\left(1-u_{[j-2]}\right)\left(\left(\partial_{j-2}-\partial_{j-1}\right)\left(1-u_{[j-2]}\right)+1\right)\E\\
	&=\left(1-u_{[j-2]}\right)\left(\partial_{j-2}-\partial_{j-1}\right)\left(1-u_{[j-2]}\right)\E +\left(1-u_{[j-2]}\right)\E .
\end{align*}
In line $3$ we used Lemma \ref{lemma1}. Let 
\begin{align*}
 \mathbb{L}^{(j)}_1&:=\left(1-u_{[j-2]}\right)\left(\partial_{j-2}-\partial_{j-1}\right) \\
 \mathbb{L}^{(j)}_2&:=\left(1-u_{[j-2]}\right)\E
\end{align*}
Both $\mathbb{L}^{(j)}_1$ and $\mathbb{L}^{(j)}_2$ can be expressed in function of the generators of $\Dm$, because of expression \eqref{vgl}.  The operator $\widetilde{C_{1j}}$ can be expressed as follows:
\[ 
	\widetilde{C_{1j}}=\mathbb{L}^{(j)}_1\mathbb{L}^{(j)}_2-(2\nu_j-1)\mathbb{L}^{(j)}_2-2\nu_1\mathbb{L}^{(j)}_1+(\nu_1+\nu_j)(\nu_1+\nu_j-1).
\]
This means that $\widetilde{C_{1j}}$ is also in $\Dm$.

Consider the first term of the operator $\widetilde{C_{2j}}$:
\begin{align*}
	&-u_{[j-2]}^2\left(1-\partial_1+\E\right)\left(\partial_{j-2}-\partial_{j-1}\right) \\
	&=-u_{[j-2]}^2\left(\partial_{j-2}-\partial_{j-1}\right)\left(-\partial_1+\E\right) \\
	&=-u_{[j-2]}\left(\left(\partial_{j-2}-\partial_{j-1}\right)u_{[j-2]}-1\right)\left(-\partial_1+\E\right) \\
	&=-u_{[j-2]}\left(\partial_{j-2}-\partial_{j-1}\right)u_{[j-2]}\left(-\partial_1+\E\right)+u_{[j-2]}\left(-\partial_1+\E\right)
\end{align*}
In line $3$ we used Lemma \ref{lemma1}. Let 
\begin{align*}
 \mathbb{L}^{(j)}_3&:=u_{[j-2]}\left(\partial_{j-2}-\partial_{j-1}\right) \\
 \mathbb{L}^{(j)}_4&:=u_{[j-2]}\left(-\partial_1+\E\right)
\end{align*}
Both $\mathbb{L}^{(j)}_3$ and $\mathbb{L}^{(j)}_4$ can be expressed in function of the generators of $\Dm$, again because of expression \eqref{vgl}.  The operator $\widetilde{C}_{2j}$ can be expressed as follows:
\[ 
	\widetilde{C}_{2j}=-\mathbb{L}^{(j)}_3\mathbb{L}^{(j)}_4-(2\nu_j-1)\mathbb{L}^{(j)}_4+2\nu_2\mathbb{L}^{(j)}_3+(\nu_2+\nu_j)(\nu_2+\nu_j-1).
\]
This means that $\widetilde{C_{2j}}$ is also in $\Dm$.

Consider the first term of the operator $\widetilde{C_{ij}}$:
\begin{align*}
&-u_{[j-1,i-2]}^2\left(\partial_{i-2}-\partial_{i-1}\right)\left(\partial_{j-2}-\partial_{j-1}\right) \\
&=-u_{[j-1,i-2]}\left(\left(\partial_{i-2}-\partial_{i-1}\right)u_{[j-1,i-2]}-1\right)\left(\partial_{j-2}-\partial_{j-1}\right) \\
&=-u_{[j-1,i-2]}\left(\partial_{i-2}-\partial_{i-1}\right)u_{[j-1,i-2]}\left(\partial_{j-2}-\partial_{j-1}\right)+u_{[j-1,i-2]}\left(\partial_{j-2}-\partial_{j-1}\right)
\end{align*}
In line $2$ we used Lemma \ref{lemma1}. Let 
\begin{align*}
 \mathbb{L}^{(ij)}_5&:=u_{[j-1,i-2]}\left(\partial_{i-2}-\partial_{i-1}\right) \\
 \mathbb{L}^{(ij)}_6&:=u_{[j-1,i-2]}\left(\partial_{j-2}-\partial_{j-1}\right)
\end{align*}
Both $\mathbb{L}^{(ij)}_5$ and $\mathbb{L}^{(ij)}_6$ can be expressed in function of the generators of $\Dm$. The operator $\widetilde{C_{ij}}$ can be expressed as follows:
\[ 
	\widetilde{C_{ij}}=-\mathbb{L}^{(ij)}_5\mathbb{L}^{(ij)}_6-(2\nu_i-1)\mathbb{L}^{(ij)}_6+2\nu_j\mathbb{L}^{(ij)}_5+(\nu_i+\nu_j)(\nu_i+\nu_j-1).
\]
This means that $\widetilde{C_{ij}}$ is also in $\Dm$.This proves Conjecture \ref{Conjecture} for this differential realization.

\section{Conclusions}

In this paper we have considered the higher rank Racah algebra in one of its differential operator realizations, obtained in the recent paper \cite{barg}. We have shown that this realization can be embedded in the enveloping algebra of a differential operator realization of $\mathfrak{sl}_n$  verifying Conjecture \ref{Conjecture} in this realization.

This embedding gives us a good idea of what should be the abstract embedding of $\mathcal{R}_{n}$ in $U(\mathfrak{sl}_{n-1})$. 
Achieving the abstract construct remains however non trivial. The differential embedding simplifies the problem since e.g. all central elements $\widetilde{C_i}$ become scalars and it is not immediate therefore how to lift the differential case to the abstract one.
An alternative construction might be better suited for that purpose and we believe that the route taken in \cite{crampe} for the case of the Heisenberg algebra is promising in this respect.

\section*{Acknowledgements}
The work of HDB is supported by the Research Foundation Flanders (FWO) under Grant EOS 30889451. HDB and WVDV are grateful for the hospitality extended to them by the Centre de Recherches Math\'ematiques in Montr\'eal, where part of this research was carried out. The research of LV is funded in part by a discovery grant of the Natural Sciences and Engineering Council (NSERC) of Canada.

%
%
%

\end{document}